\documentclass[graybox]{svmult} % this had errors for %equations
\usepackage{mathptmx}        % selects Times Roman as basic font
\usepackage{helvet}          % selects Helvetica as sans-serif font
\usepackage{courier}         % selects Courier as typewriter font
\usepackage{makeidx}         % allows index generation
\usepackage{graphicx}        % standard LaTeX graphics tool
\usepackage{multicol}        % used for the two-column index
\usepackage[bottom]{footmisc}% places footnotes at page bottom

\usepackage{color}

\definecolor{hhbblue}{rgb}{.8, .8, 1}
  \newcommand*\hhbbluebox[1]{%
    \colorbox{hhbblue}{\hspace{1em}#1\hspace{1em}}}

\usepackage{soul} % HHB : added to communicate to authors
\usepackage{empheq}
\usepackage{enumerate}
\usepackage{enumitem}

\usepackage{placeins}

\makeindex             % used for the subject index
\usepackage{amsmath}
\usepackage{amssymb}
\usepackage{texdraw}
\usepackage{bbm}
\usepackage{color}
\usepackage{etex}
\usepackage{tikz,pgfplots}
\usepackage{booktabs}

\usepackage{subfigure}
\usepackage{exscale,relsize}
\usepackage{pifont}
\usepackage{psfrag}
\usepackage{amsbsy}
\usepackage{latexsym}

\spnewtheorem{thm}[theorem]{Theorem}{\bfseries}{\rmfamily}
\spnewtheorem{lem}[theorem]{Lemma}{\bfseries}{\rmfamily}
\spnewtheorem{cor}[theorem]{Corollary}{\bfseries}{\rmfamily}
\spnewtheorem{propn}[theorem]{Proposition}{\bfseries}{\rmfamily}
\spnewtheorem{defn}[theorem]{Definition}{\bfseries}{\rmfamily}
\spnewtheorem{rem}[theorem]{Remark}{\bfseries}{\rmfamily}
\spnewtheorem{ex}[theorem]{Example}{\bfseries}{\rmfamily}

\spdefaulttheorem{hhbtheorem}{Theorem}{\bfseries}{\rmfamily}% [section]
\spnewtheorem{hhblemma}[hhbtheorem]{Lemma}{\bfseries}{\rmfamily}
\spnewtheorem{hhbcorollary}[hhbtheorem]{Corollary}{\bfseries}{\rmfamily}
\spnewtheorem{hhbproposition}[hhbtheorem]{Proposition}{\bfseries}{\rmfamily}
\spnewtheorem{hhbdefinition}[hhbtheorem]{Definition}{\bfseries}{\rmfamily}
\spnewtheorem{hhbquestion}[hhbtheorem]{Open Question}{\bfseries}{\rmfamily}
\spnewtheorem{hhbassumption}[hhbtheorem]{Assumption}{\bfseries}{\rmfamily}%}
\spnewtheorem{hhbcondition}[hhbtheorem]{Condition}{\bfseries}{\rmfamily}%}
\spnewtheorem{hhbalgorithm}[hhbtheorem]{Algorithm}{\bfseries}{\rmfamily}%}
\spnewtheorem{hhbexample}[hhbtheorem]{Example}{\bfseries}{\rmfamily}%}
\spnewtheorem{hhbfact}[hhbtheorem]{Fact}{\bfseries}{\rmfamily}
\spnewtheorem{hhbremark}[hhbtheorem]{Remark}{\bfseries}{\rmfamily}%}

\newcommand{\hhbscal}[2]{\langle{{#1},{#2}}\rangle}
\newcommand{\hhbHH}{\ensuremath{X}}
\newcommand{\hhbTo}{\ensuremath{\rightrightarrows}}
\newcommand{\hhbgr}{\ensuremath{\operatorname{gra}}}
\newcommand{\hhbId}{\ensuremath{\operatorname{Id}}}
\newcommand{\hhbzer}{\ensuremath{\operatorname{zer}}}
\newcommand{\hhbweakly}{\ensuremath{\,\rightharpoonup}\,}
\newcommand{\hhbmenge}[2]{\big\{{#1} \mid {#2}\big\}}
\newcommand{\hhbNN}{\ensuremath{\mathbb N}}
\newcommand{\hhbnnn}{\ensuremath{{n \in \hhbNN}}}
\newcommand{\hhbthalb}{\ensuremath{\tfrac{1}{2}}}
\newcommand{\hhbbX}{\ensuremath{\mathbf{X}}}
\newcommand{\hhbbC}{\ensuremath{\mathbf{C}}}
\newcommand{\hhbbD}{\ensuremath{\mathbf{D}}}
\newcommand{\hhbbF}{\ensuremath{\mathbf{F}}}
\newcommand{\hhbby}{\ensuremath{\mathbf{y}}}
\newcommand{\hhbbx}{\ensuremath{\mathbf{x}}}
\newcommand{\hhbbz}{\ensuremath{\mathbf{z}}}
\newcommand{\hhbFix}{\ensuremath{\operatorname{Fix}}}
\newcommand{\hhbcard}{\ensuremath{\operatorname{card}}}

\renewcommand\theenumi{\rm (\roman{enumi})}

\begin{document}

\title{The method of cyclic intrepid projections:\\ convergence analysis and
numerical experiments}

% following used only if title too long
\titlerunning{The method of cyclic intrepid projections} 

\author{Heinz H.\ Bauschke 
\and 
Francesco Iorio
\and 
Valentin R.\ Koch 
}
\institute{
Heinz H.\ Bauschke 
\at 
Mathematics, 
University of British Columbia, 
Kelowna, BC~ V1V 1V7, Canada,
\email{heinz.bauschke@ubc.ca}
\and
Francesco Iorio
\at 
Autodesk Research,
210 King Street East,
Suite 600,
Toronto, ON~ M5A 1J7, Canada,
\email{francesco.iorio@autodesk.com}
\and 
Valentin R.\ Koch
\at 
Information Modeling \& Platform Products Group (IPG),
Autodesk, Inc.,
\email{valentin.koch@autodesk.com}
}

\maketitle

\medskip

% duplicate abstract in both.
\abstract{
The convex feasibility problem asks to find a point in the
intersection of a collection of nonempty closed convex sets. 
This problem is of basic importance in mathematics and the
physical sciences, and projection
(or splitting) methods solve it by employing the projection
operators associated with the individual sets to generate a
sequence which converges to a solution. 
Motivated by an application in road design, we present
the method of cyclic intrepid projections (CycIP) and provide 
a rigorous convergence analysis. We also report on very promising 
numerical experiments in which CycIP is compared to a commerical
state-of-the-art optimization solver. 
}

\abstract*{
The convex feasibility problem asks to find a point in the
intersection of a collection of nonempty closed convex sets. 
This problem is of basic importance in mathematics and the
physical sciences, and projection
(or splitting) methods solve it by employing the projection
operators associated with the individual sets to generate a
sequence which converges to a solution. 
Motivated by an application in road design, we present
the method of cyclic intrepid projections (CycIP) and provide 
a rigorous convergence analysis. We also report on very promising 
numerical experiments in which CycIP is compared to a commerical
state-of-the-art optimization solver. 
}



\section{Introduction}

\label{hhbsec:Intro}

\noindent
Throughout this paper, we assume that
\begin{empheq}[box=\hhbbluebox]{equation}
\text{ $X$ is a real Hilbert space} 
\end{empheq}
with inner product $\hhbscal{\cdot}{\cdot}$ and induced norm $\|\cdot\|$.
(We also write $\|\cdot\|_2$ instead of $\|\cdot\|$ if we
wish to emphasize this norm compared to other norms.
We assume basic notation and results from
convex analysis and fixed point theory; see, e.g., 
\cite{hhbBC,hhbBorVan,hhbGoeKir,hhbGoeRei,hhbRockWets}.) 

Let $(C_i)_{i\in I}$ be a finite family of closed convex subsets
of $X$ such that
\begin{equation}
C := \bigcap_{i\in I} C_i\neq\varnothing. 
\end{equation}
We aim to find a point in $C$ given that the individual
constraint sets $C_i$ are simple in the sense that their
associated projections\footnote{
Given a nonempty subset $S$ of $X$ and $x\in X$, we write $d_S(x) :=
\inf_{s\in S}\|x-s\|$ for the \emph{distance} from $x$ to $S$.
If $S$ is also closed and convex, then the infimum defining
$d_S(x)$ is attained at a \emph{unique} vector called the 
\emph{projection} of $x$ onto $S$ and denoted by $P_S(x)$ or
$P_Sx$.} are easy to compute. To solve the widespread 
\emph{convex feasibility problem}
``find $x\in C$'', we employ \emph{projection methods}.
These splitting-type methods use the individual projections
$P_{C_i}$ in order to generate a sequence that converges to a
point in $C$. For further information, we refer the
reader to, e.g., 
\cite{hhbSIREV,hhbBC,hhbCegielski,hhbCenZen,hhbPLC97,hhbGoeKir,hhbGoeRei,hhbHerman09}. 

In previous work on a feasibility problem arising in road
design \cite{hhbBK}, the \emph{method of cyclic intrepid
projections (CycIP)} was found to be an excellent overall algorithm. 
Unfortunately, CycIP was applied heuristically without an
underlying convergence result. 

{\em 
The goal of this paper is two-fold. First, 
we present a checkable condition sufficient for convergence and
provide a rigorous convergence proof. In fact, our main result
applies to very general feasibility problems satisfying an
interiority assumption. Second, we numerically compare CycIP
to a commercial LP solver for test problems that are both convex
and nonconvex to evaluate competitiveness of CycIP. 
}

The remainder of the paper is organized as follows. 
In Section~\ref{hhbsec:Prelim} we provide basic properties of
projection operators. 
Useful results on Fej\'er monotone sequences are recalled in
Section~\ref{hhbsec:Fejer}. 
Our main convergence results are presented in
Section~\ref{hhbsec:main}. 
In Section~\ref{hhbsec:road}, we review the feasibility
problem arising in road design and obtain a rigorous convergence
result for CycIP. 
We report on numerical experiments in 
Section~\ref{hhbsec:gurobi} and offer concluding remarks in 
Section~\ref{hhbsec:conclusion}. 

We end this section with notation.
The closed ball centered at $y\in X$ of radius $r$ is
$B(y;r) := \hhbmenge{z\in X}{\|y-z\|\leq r}$.
Finally, we write $\mathbb{R}_+$ and $\mathbb{R}_{++}$ 
for the nonnegative real numbers and strictly positive reals,
respectively.

\section{Relaxed and intrepid projectors}

\label{hhbsec:Prelim}

In this section, we introduce the key operators used in the
projection methods studied later.

\begin{hhbfact}[relaxed projector]
\label{hhbf:relax}
Let $C$ be a nonempty closed convex subset of $X$,
and let $\lambda \in \left]0,2\right[$. 
Set $R := (1-\lambda)\hhbId + \lambda P_C$,
let $x\in X$, and let $c\in C$.
Then
\begin{equation}
\|x-c\|^2 - \|Rx-c\|^2 \geq
\frac{2-\lambda}{\lambda}\|x-Rx\|^2
= (2-\lambda)\lambda d_C^2(x).
\end{equation}
\end{hhbfact}
\begin{proof}
Combine \cite[Lemma~2.4.(iv)]{hhbSIREV} with 
\cite[Proposition~4.8]{hhbBC}. 
\qed
\end{proof}

In fact, the relaxed projector is an example of a so-called
\emph{averaged map}; see, e.g., \cite{hhbBBR,hhbBC,hhbCombettes04} 
for more on this useful notion.

\begin{hhbdefinition}[enlargement]
Given a nonempty closed convex subset $Z$ of $X$, and
$\alpha\in\mathbb{R}_+ := \hhbmenge{\xi\in\mathbb{R}}{\xi \geq
0}$, we write
\begin{equation}
C_{[\alpha]} := \hhbmenge{x\in X}{d_C(x)\leq \alpha}  = C +
B(0;\alpha)
\end{equation}
and call $C_{[\alpha]}$ the \emph{$\alpha$-enlargement} of $C$. 
\end{hhbdefinition}

Note that $C_{[0]} = C$,
that $C_{[\alpha]}$ is a nonempty closed convex subset of $X$,
and that
if $\alpha<\beta$, then
$C_{[\alpha]}\subseteq C_{[\beta]}$. 
We mention in passing that the \emph{depth}\footnote{This
function is considered, e.g., in \cite[Exercise~8.5]{hhbBV}.} 
of each $z\in Z$ (with respect to $C$), 
i.e., $d_{X\smallsetminus C}(z)$, is at least $\alpha$.

\begin{hhbfact}
\label{hhbf:blowproj}
{\rm (See, e.g., \cite[Proposition~28.10]{hhbBC}.)}
Let $C$ be a nonempty closed convex subset of $X$,
and let $\beta\in\mathbb{R}_{+}$. 
Set $D := C_{[\beta]}$. Then
\begin{equation}
(\forall x\in X)\quad
P_D x = 
\begin{cases}
x, &\text{if $d_C(x)\leq \beta$;}\\
\displaystyle P_Cx + \beta\frac{x-P_Cx}{d_C(x)}, &\text{otherwise.}
\end{cases}
\end{equation}
\end{hhbfact}

\begin{hhbdefinition}[intrepid projector]
Let $Z$ be a nonempty closed convex subset of $X$,
let $\beta\in\mathbb{R}_+$, and set
$C := Z_{[\beta]}$.
The corresponding 
\emph{intrepid projector $Q := Q_C$ onto $C$ (with respect to $Z$ and
$\beta$)} is defined by 
\begin{align}
Q\colon X\to X\colon
x &\mapsto 
x+
\Big(1-\tfrac{P_{[\beta,2\beta]}d_Z(x)}{\beta}\Big)\big(x-P_Zx\big)
\notag\\
&= 
\begin{cases}
P_Zx, &\text{if $d_Z(x)\geq 2\beta$;}\\
x, &\text{if $d_Z(x) \leq \beta$;}\\
\displaystyle x+ \Big(1-\tfrac{d_Z(x)}{\beta}\Big)\big(x-P_Zx\big), &\text{otherwise.}
\end{cases} 
\end{align}
We refer to these three steps as the
\emph{projection step},
the \emph{identity step}, and
the \emph{reflection step}, respectively. 
\end{hhbdefinition}

\begin{hhbexample}[intrepid projector onto a hyperslab \`a la
Herman]
Suppose that $a\in X\smallsetminus\{0\}$,
let $\alpha\in\mathbb{R}$, 
let $\beta\in\mathbb{R}_+$, 
and set $Z := \hhbmenge{x\in X}{\hhbscal{a}{x}=\alpha}$.
Then $Z$ is a hyperplane and $Z_{[\beta]}$ is a hyperslab. 
Moreover, the associated intrepid projector onto $Z_{[\beta]}$ is precisely
the operator considered by Herman in \cite{hhbHerman75}.
\end{hhbexample}

\begin{hhbproposition}[basic properties of the intrepid
projector]
\label{hhbp:bpip}
Let $Z$ be a nonempty closed convex subset of $X$,
let $\beta\in\mathbb{R}_+$, set
$C := Z_{[\beta]}$, and
denote the corresponding 
intrepid projector onto $C$ (with respect to $Z$ and
$\beta$) by $Q$.
Now 
let $\alpha\in[0,\beta]$,
and let $y \in Z_{[\alpha]}$, and
let $x\in X$. 
Then 
\begin{equation}
Qx \in [x,P_Zx]\cap C
\end{equation}
and 
exactly one of the following holds: 
\begin{enumerate}
\item 
\label{hhbp:bpip1}
$d_Z(x)\leq \beta$, $x=Qx\in C$, and 
$\|x-y\|^2 -\|Qx-y\|^2 = 0$. 
\item 
\label{hhbp:bpip2}
$d_Z(x)\geq 2\beta$ and
\begin{align}
\|x-y\|^2 - \|Qx-y\|^2 &\geq 2(\beta-\alpha)\|x-Qx\| =
2(\beta-\alpha)d_Z(x) \notag\\
&=2(\beta-\alpha)\big(\beta+d_C(x)\big)\geq 4\beta(\beta-\alpha).
\end{align}
\item 
\label{hhbp:bpip3}
$\beta < d_Z(x) < 2\beta$ and 
\begin{align}
\|x-y\|^2 - \|Qx-y\|^2 &\geq 2(\beta-\alpha)\|x-Qx\|
= \frac{2(\beta-\alpha)}{\beta}d_Z(x)\big(d_Z(x)-\beta\big)\notag\\
&=
\frac{2(\beta-\alpha)}{\beta}d_C(x)\big(\beta+d_C(x)\big).
\end{align}
\end{enumerate}
Consequently, in \emph{every} case, we have
\begin{equation}
\|x-y\|^2 - \|Qx-y\|^2 \geq 2(\beta-\alpha)\|x-Qx\|
\geq 2(\beta-\alpha)d_C(x). 
\end{equation}
\end{hhbproposition}
\begin{proof}
Set $\delta := d_Z(x)$, $p := P_Zx$, and 
write $y= z+\alpha b$, where $z\in Z$,  $b\in X$ and $\|b\|\leq 1$. 
Note that if $\delta > \beta$, then $d_C(x)=\delta-\beta$ using 
Fact~\ref{hhbf:blowproj}. 

\ref{hhbp:bpip1}:
This follows immediately from the definition of $Q$.

\ref{hhbp:bpip2}:
Using Cauchy--Schwarz in \eqref{hhbe:1123a}, 
and the projection theorem (see, e.g.,
\cite[Theorem~3.14]{hhbBC}) in \eqref{hhbe:1123b},
we obtain
{\allowdisplaybreaks
\begin{subequations}
\begin{align}
\|x-y\|^2 - \|Qx-y\|^2 &=
\|x-(z+\alpha b)\|^2 - \|p-(z+\alpha b)\|^2 \notag \\
&= \|x\|^2 - \|p\|^2 -2\hhbscal{x}{z+\alpha b} +
2\hhbscal{p}{z+\alpha b}\notag \\
&\geq \|x\|^2 - \|p\|^2 + 2\hhbscal{p-x}{z} 
-2\alpha\|p-x\|\|b\|\label{hhbe:1123a}\\
&\geq \|x\|^2-\|p\|^2 + 2\hhbscal{p-x}{z-p} + 2\hhbscal{p-x}{p} -
2\alpha\delta\notag\\
&\geq \|x\|^2-\|p\|^2  + 2\hhbscal{p-x}{p} - 
2\alpha\delta\label{hhbe:1123b}\\
&= \|x\|^2 + \|p\|^2 - 2\hhbscal{x}{p} - 
2\alpha\delta
= \|x-p\|^2 - 2\alpha\delta\notag\\
&=\delta^2-2\alpha\delta
=\delta\big(\delta-2\alpha\big)
\geq 2\delta(\beta-\alpha) 
\geq 4\beta(\beta-\alpha). \notag
\end{align}
\end{subequations}
}

\ref{hhbp:bpip3}:
Set $\eta := (\delta-\beta)/\beta\in\left]0,1\right[$.
Then $Qx = (1-\eta)x + \eta p$ and hence
\begin{equation}
\label{hhbe:1123d}
\|x-Qx\| = \eta\|x-p\|=\eta\delta =
\frac{\delta-\beta}{\beta}\delta.
\end{equation}
Using, e.g.,\ 
\cite[Corollary~2.14]{hhbBC} in \eqref{hhbe:1123c1}, 
and Cauchy--Schwarz 
and \cite[Theorem~3.14]{hhbBC} in \eqref{hhbe:1123c2},
we obtain 
{\allowdisplaybreaks
\begin{subequations}
\label{hhbe:1123c}
\begin{align}
\|x-y\|^2 - \|Qx-y\|^2 &=
\|x\|^2 - \|Qx\|^2 -2\hhbscal{x}{y} +
2\hhbscal{Qx}{y}\notag \\
&=\|x\|^2 - \| (1-\eta)x + \eta p\|^2 
+ 2\hhbscal{(1-\eta)x+\eta p - x}{y}\notag\\
& = \|x\|^2 - (1-\eta)\|x\|^2 - \eta\|p\|^2
\label{hhbe:1123c1}\\
&\qquad
+\eta(1-\eta)\|x-p\|^2
+ 2\eta\hhbscal{p-x}{z+\alpha b} \notag\\
&=\eta\big(\|x\|^2-\|p\|^2\big) + (1-\eta)\eta\|x-p\|^2 \notag\\
&\qquad + 2\eta \big( \hhbscal{p-x}{z-p} + \hhbscal{p-x}{p} +
\alpha\hhbscal{p-x}{b}\notag\big)\\
&\geq \eta\big(\|x\|^2-\|p\|^2 +(1-\eta)\|x-p\|^2 +
2\hhbscal{p-x}{p}\big)\label{hhbe:1123c2}\\
&\qquad -2\alpha\eta\|x-p\|\notag\\
&=\eta\big( \|x-p\|^2 + (1-\eta)\|x-p\|^2
-2\alpha\|x-p\|\big)\notag\\
&=\delta\eta\big((2-\eta)\delta -2\alpha\big)\notag\\
&=\frac{\delta}{\beta}(\delta-\beta)\Big(\big(2-(\delta-\beta)\beta^{-1}\big)\delta
- 2\alpha\Big)\notag\\
&=\frac{\delta}{\beta}(\delta-\beta)\big( -\beta^{-1}\delta^2 +
3\delta -2\alpha\big).\notag
\end{align}
\end{subequations}
}
Now the quadratic
$q\colon [\beta,2\beta]\to\mathbb{R}\colon
\xi \mapsto  -\beta^{-1}\xi^2 +
3\xi-2\alpha$ has a maximizer at $\xi = (3/2)\beta$ and
it satisfies $q(\beta) = q(2\beta) = 2(\beta-\alpha) \geq 0$. 
It follows that $\min q\big([\beta,2\beta]\big) =
2(\beta-\alpha)$. 
Therefore, \eqref{hhbe:1123c} and \eqref{hhbe:1123d} 
\begin{equation}
\|x-y\|^2 - \|Qx-y\|^2 
\geq 
\frac{\delta}{\beta}(\delta-\beta)2(\beta-\alpha)
= 2(\beta-\alpha)\|x-Qx\|. 
\end{equation}
The proof of the ``Consequently'' part follows easily. 
\qed
\end{proof}

Proposition~\ref{hhbp:bpip} implies that
the intrepid projector is \emph{quasi nonexpansive}; see, e.g.,
\cite{hhbMOR,hhbCegielski,hhbSY,hhbYO,hhbYYY} for further results utilizing this notion. 

\section{Fej\'er monotonicity}

\label{hhbsec:Fejer}

We now review the definition and basic results on
Fej\'er monotone sequences. These will be useful in establishing
our convergence results. 

\begin{hhbdefinition}
Let $(x_k)_{k\in\mathbb{N}}$ be a sequence in $X$,
and let $C$ be a nonempty closed convex subset of $X$.
Then  $(x_k)_{n\in\mathbb{N}}$ is
\emph{Fej\'er monotone with respect to $C$} if 
\begin{equation}
(\forall k\in\mathbb{N})(\forall c\in C)\quad
\|x_{k+1}-c\|\leq\|x_{k}-c\|.
\end{equation}
\end{hhbdefinition}

\begin{hhbfact}
\label{hhbf:Fejer}
Let $(x_k)_{k\in\mathbb{N}}$ be a sequence in $X$ that
is Fej\'er monotone with respect to some nonempty closed convex
subset $C$ of $X$. Then the following hold:
\begin{enumerate}
\item 
\label{hhbf:Fejer1}
If $\operatorname{int} C\neq\varnothing$, then
$(x_k)_{k\in\mathbb{N}}$ converges strongly to some point in $X$.
\item 
\label{hhbf:Fejer2}
If each weak cluster point of  $(x_k)_{k\in\mathbb{N}}$ lies in
$C$, then $(x_k)_{k\in\mathbb{N}}$
converges weakly to some point in $C$. 
\item 
\label{hhbf:Fejer3}
If $d_C(x_k)\to 0$, then
$(x_k)_{k\in\mathbb{N}}$ converges strongly to some point in $C$. 
\end{enumerate}
\end{hhbfact}

\begin{proof}
See, e.g., 
\cite{hhbSIREV}, \cite[Chapter~5]{hhbBC}, or \cite{hhbCombettes01}. 
\qed
\end{proof}

\section{The method of cyclic intrepid projections}

\label{hhbsec:main}

We now assume that each $C_i$ is a closed convex subset
of $X$, with 
\begin{empheq}[box=\hhbbluebox]{equation}
C := \bigcap_{i\in I} C_i\neq\varnothing.
\end{empheq}
The index set is split into two sets, corresponding
to enlargements and regular sets: 
\begin{empheq}[box=\hhbbluebox]{equation}
I_0 := \hhbmenge{i\in I}{C_i = (Z_i)_{[\beta_i]},\;\text{where
$\beta_i>0$}}
\quad\text{and}\quad
I_1 := I\smallsetminus I_0.
\end{empheq}
Assume an index selector map
\begin{empheq}[box=\hhbbluebox]{equation}
\mathrm{i}\colon \mathbb{N}\to I,
\end{empheq}
where $(\forall i\in I)$ $\mathrm{i}^{-1}(i)$ is an infinite
subset of $\mathbb{N}$. 
We say that the \emph{control is quasicyclic with quasiperiod $M
\in\{1,2,\ldots\}$}
if 
$(\forall k\in\mathbb{N})$
$I = \{
\mathrm{i}(k),\mathrm{i}(k+1),\ldots,\mathrm{i}(k+M-1)\}$. 

Let $(\lambda_i)_{i\in I_1}$ be a family in $\left]0,2\right[$. 
We define a family of operators 
\begin{empheq}[box=\hhbbluebox]{equation}
(T_i)_{i\in I}
\end{empheq}
from $X$ to $X$ as follows.
If $i\in I_0$, then $T_i$ is the intrepid projector onto $C_i$
(with respect to $Z_i$ and $\beta_i$);
if $i\in I_1$, then $T_i$ is the relaxed projector onto $C_i$
with relaxation parameter $\lambda_i$. 

\begin{hhbalgorithm}[method of cyclic intrepid projections]
\label{hhba:cycip}
Given a starting point $x_0\in X$, 
the \emph{method of intrepid projections} proceeds via
\begin{empheq}[box=\hhbbluebox]{equation}
(\forall k\in\mathbb{N})\quad
x_{k+1} := T_{\mathrm{i}(k)} x_k.
\end{empheq}
\end{hhbalgorithm}

We begin our analysis with a simple yet useful observation. 

\begin{hhblemma}
\label{hhbl:Fejer}
The sequence $(x_k)_{k\in\mathbb{N}}$ generated by
Algorithm~\ref{hhba:cycip} is
Fej\'er monotone with respect to $C$.
\end{hhblemma}
\begin{proof}
Combine Fact~\ref{hhbf:relax} with 
Proposition~\ref{hhbp:bpip}. 
\qed
\end{proof}

We now deepen our convergence analysis. 
We start with the purely intrepid case.

\begin{hhbtheorem}[intrepid projections only]
\label{hhbt:ipo}
Suppose that $I_1=\varnothing$ and that 
$\operatorname{int} C\neq\varnothing$. Then
$(x_k)_{k\in\mathbb{N}}$ converges strongly to some point in $C$.
\end{hhbtheorem}

\begin{proof}
Combining Lemma~\ref{hhbl:Fejer} 
with Fact~\ref{hhbf:Fejer}\ref{hhbf:Fejer1}, 
we deduce that 
$(x_k)_{k\in\mathbb{N}}$ converges strongly to some point $\bar{x}\in
X$. 
Let $i\in I$. Then
the subsequence $(x_{\mathrm{i}^{-1}(i)})$ not only lies in $C_i$ (see
Proposition~\ref{hhbp:bpip}) but it also converges to $\bar{x}$.
Since $C_i$ is closed, we deduce that $\bar{x}\in C_i$.  
\qed
\end{proof}

The proof of the following result follows that of Herman \cite{hhbHerman75}
who considered more restrictive controls. 
(See also \cite{hhbChenHerman}.)

\begin{hhbcorollary}[parallelotope]
\label{hhbc:Herman}
Suppose that $X$ is finite-dimensional, 
that $I_1=\varnothing$, 
that each $Z_i$ is a hyperplane, and that
$\operatorname{int} C\neq\varnothing$. 
Then $(x_k)_{k\in\mathbb{N}}$ converges to some point in 
the parallelotope $C$ \emph{in finitely many steps}. 
\end{hhbcorollary}

\begin{proof}
By Theorem~\ref{hhbt:ipo}, $(x_k)_{k\in\mathbb{N}}$ converges
to some point $\bar{x}\in C$. 
If $\hhbmenge{x_k}{k\in\mathbb{N}} \cap \operatorname{int}
C\neq\varnothing$, then  $(x_k)_{k\in\mathbb{N}}$ is eventually
constant. 
Assume to the contrary that  $(x_l)_{k\in\mathbb{N}}$ is not
eventually constant. 
Then $\bar{x}\notin \hhbmenge{x_k}{k\in\mathbb{N}}\cup 
\operatorname{int} C$.
Since each $C_i$ is a hyperslab,
$ \operatorname{bdry} C_i$ is the union of two disjoint
hyperplanes parallel to $Z_i$. 
We collect these finitely many hyperplanes in a set $H$.
The finite collection of these hyperplanes containing $\bar{x}$,
which we denote by 
$H(\bar{x})$, is nonempty. 
Moreover, $(x_k)_{k\in\mathbb{N}}$ cannot have arisen with
infinitely many projection steps as these only occur at a minimum
distance from the sets. 
Therefore, infinitely many reflection steps have been executed. 
Hence there exists $K_1\in\mathbb{N}$ such that iteration index
$k$ onwards, we only execute identity or reflection steps. 
Now let $\varepsilon>0$ be sufficiently small 
such that $B(\bar{x};\varepsilon)$
makes an empty intersection with every hyperplane drawn from
$H\smallsetminus H(\bar{x})$. 
Since $x_k\to\bar{x}$, there exists $K_2\in\mathbb{N}$ such that
$(\forall k\geq N_2)$ $x_k\in B(\bar{x};\varepsilon)$. 
Since $\bar{x}\in C$ and $(x_k)_{k\in\mathbb{N}}$ is Fej\'er
monotone with respect to $C$, it follows that
$(\forall k\in\mathbb{N})$
$\|x_{k+1}-\bar{x}\|\leq\|x_k-\bar{x}\|$. 
Hence the aforementioned reflection steps from $K_2$ onwards
must be all with respect to hyperplanes taken from $H(\bar{x})$.
Set $K := \max\{K_1,K_2\}$.
It follows altogether that
$(\forall k\geq K)$
$0<\|x_{k+1}-\bar{x}\|=\|x_k-\bar{x}\|$.
But this is absurd since $x_k\to \bar{x}$.
\qed 
\end{proof}

\begin{hhbremark}
Both Theorem~\ref{hhbt:ipo}
and Corollary~\ref{hhbc:Herman}
fail if $\operatorname{int} C=\varnothing$:
indeed, consider two hyperslabs $C_1$ and $C_2$ in $\mathbb{R}^2$
such that $C$ is a line (necessarily parallel to $Z_1$ and $Z_2$). 
If we start the iteration
sufficiently close to this line, but not on this line, then
$(x_n)_{n\in\mathbb{N}}$ will oscillate between two point outside
$C$. 

Furthermore, finite convergence may fail in Corollary~\ref{hhbc:Herman}
without the interiority assumption: indeed, consider a hyperslab
in $\mathbb{R}^2$ which is intersected by a line at an angle
strictly between $0$ and $\pi/4$. 
\end{hhbremark}

We now present our fundamental convergence result.

\begin{hhbtheorem}[main result]
\label{hhbt:main}
Suppose that 
$\bigcap_{i\in I_1} C_i \cap \bigcap_{i\in I_0}\operatorname{int}C_i
\neq \varnothing$ and that the control is
quasicyclic. 
Then the sequence $(x_k)_{k\in\mathbb{N}}$ generated by
Algorithm~\ref{hhba:cycip} converges weakly to some point in
$C$. 
The convergence is strong provided one of the following
conditions holds:
\begin{enumerate}
\item 
\label{hhbt:main1}
$X$ is finite-dimensional.
\item
\label{hhbt:main2}
$I_1$ is either empty or a singleton.
\end{enumerate}
\end{hhbtheorem}

\begin{proof}
By Lemma~\ref{hhbl:Fejer}, $(x_k)_{k\in\mathbb{N}}$ is 
Fej\'er monotone with respect to $C$. 
Take $y\in \bigcap_{i\in I_1} C_i 
\cap \bigcap_{i\in I_0}\operatorname{int}C_i $.
Writing $\|x_0-y\|^2 =
\sum_{k\in\mathbb{N}}\|x_k-y\|^2-\|x_{k+1}-y\|^2$,
and recalling Fact~\ref{hhbf:relax} and
Proposition~\ref{hhbp:bpip}, we deduce that
$x_k-x_{k+1}\to 0$ and that 
$d_{C_{\mathrm{i}(n)}}(x_k)\to 0$.
The quasicyclicality of the control now yields
\begin{equation}
\max\hhbmenge{d_{C_i}(x_k)}{i\in I}\to 0.
\end{equation}
Therefore, every weak cluster point of $(x_k)_{k\in\mathbb{N}}$
lies in $C$. 
By Fact~\ref{hhbf:Fejer}\ref{hhbf:Fejer2}, 
there exists $\bar{x}\in X$ such that
\begin{equation}
x_k\rightharpoonup \bar{x} \in C
\end{equation}
as announced. 

Let us turn to strong convergence. 
Item~\ref{hhbt:main1} is obvious since strong and weak
convergence coincide in finite-dimensional Hilbert space.

Now consider \ref{hhbt:main2}.
If $I_1=\varnothing$, then strong convergence follows from 
Theorem~\ref{hhbt:ipo}.
Thus assume that $I_1$ is a singleton. 
By, e.g., \cite[Theorem~5.14]{hhbSIREV}, 
\begin{equation}
d_C(x_k)\to 0.
\end{equation}
Hence, using Fact~\ref{hhbf:Fejer}\ref{hhbf:Fejer3}, 
we conclude that $x_k\to \bar{x}$. 
\qed
\end{proof}

\begin{hhbremark}
Our sufficient conditions for strong convergence are sharp:
indeed, Hundal's example \cite{hhbHundal} shows that
strong convergence may fail if
{\rm (i)} $X$ is infinite-dimensional and
{\rm (ii)} $I_1$ contains more than one element. 
\end{hhbremark} 

\section{CycIP and the road design problem}

\label{hhbsec:road}

From now on,
we assume that 
\begin{empheq}[box=\hhbbluebox]{equation}
\text{ $X=\mathbb{R}^n$,} 
\end{empheq}
and that we are given $n$ breakpoints
\begin{empheq}[box=\hhbbluebox]{equation}
t = (t_1,\ldots,t_n)\in X
\quad\text{such that~} 
t_1 < \cdots < t_n.
\end{empheq}
The problem is to 
\begin{empheq}[box=\hhbbluebox]{equation}
\text{find $x=(x_1,\ldots,x_n) \in X$}
\end{empheq}
such that all of the following constraints are satisfied:
\begin{itemize}
\item
\textbf{interpolation constraints:}
For a subset $J$ of $\{1,\ldots,n\}$, we have $x_j=y_j$, where
$y\in\mathbb{R}^J$ is given. 
\item
\textbf{slope constraints:}
each slope $s_j := 
(x_{j+1}-x_j)/(t_{j+1}-t_j)$ satisfies
$|s_j|\leq \sigma_j$, where 
$j\in\{1,\ldots,n-1\}$ and 
$\sigma \in \mathbb{R}^{n-1}_{++}$ is given. 
\item
\textbf{curvature constraints:}
$\gamma_j \geq s_{j+1}-s_j \geq \delta_j$,
for every $j\in\{1,\ldots,n-2\}$, and for given 
$\gamma$ and $\delta$  in $\mathbb{R}^{n-2}$.
\end{itemize}
This problem is of fundamental interest in road design; see \cite{hhbBK}
for further details. 

By grouping the constraints appropriately, this feasibility problem can be
reformulated as the following convex feasibility problem involving six sets:
\begin{empheq}[box=\hhbbluebox]{equation}
\label{hhbe:road}
\text{find $x\in C = \bigcap_{i\in I} C_i = C_1\cap \cdots \cap C_6$,}
\end{empheq}
where $I:=\{1,\ldots,6\}$; see \cite[Section~2]{hhbBK} for details. 
These sets have additional structure: $C_1$ is an affine subspace
incorporating the interpolation constraints,
$C_2$ and $C_3$ are both intersections of hyperslabs with normal vectors
having disjoint support modeling the slope constraints, and the
curvature constraints are similarly incorporated through 
$C_4$, $C_5$, and $C_6$. 
All these sets have explicit and easy-to-implement (regular and intrepid)
projection formulas. 
Since only the set $C_1$ has no interior, we set $T_1 = P_{C_1}$.
For every $i\in\{2,\ldots,6\}$, we 
set $Q_{C_i}$.
(If we set $T_i = P_{C_i}$, we get the classical method of cyclic
projections.)
This gives rise to the algorithm, which we call the method of 
\textbf{cyclic intrepid projections (CycIP)}. 

We thus obtain the following consequence of our main result
(Theorem~\ref{hhbt:main}): 

\begin{hhbcorollary}[strict feasibility]
\label{hhbc:road}
Suppose that 
$C_1 \cap \bigcap_{2\leq i\leq 6}\operatorname{int} C_i\neq \varnothing$, i.e.,
there exists a \emph{strictly feasible} solution to \eqref{hhbe:road},
i.e., it satisfies 
the interpolation constraints, and it satisfies the 
slope and curvature constraint inequalities \emph{strictly}. 
Then the sequence generated by CycIP converges to
a solution of \eqref{hhbe:road}. 
\end{hhbcorollary}

In \cite{hhbBK}, which contains a comprehensive comparison of various
algorithms for solving \eqref{hhbe:road}, CycIP was found to be the best
overall algorithm. However, due to the interpolation constraint set $C_1$,
which has \emph{empty interior}, the convergence of CycIP is not 
guaranteed by Theorem~\ref{hhbt:ipo} or convergence results
derived earlier. 
Corollary~\ref{hhbc:road} is the first \emph{rigorous
justification} of CycIP in the setting of road design.

\begin{hhbremark}[nonconvex minimum-slope constraints]
\label{hhbrem:nonconvex}
In \cite{hhbBK} we also considered a variant of the slope constraints with
an imposed minimal strictly positive slope. This is a setting of
significant interest in road design as zero slopes are not favoured because
of, e.g., drainage problems. The accordingly 
modified sets $C_2$ and $C_3$ are in that case
\emph{nonconvex}; however, explicit formulas for (regular and
intrepid) projections are still available. 
The application of CycIP must
then be regarded as a heuristic as there is no accompanying body of
convergence results.
\end{hhbremark}

In the following section, we will investigate the numerical
performance of CycIP and compare it to a linear programming solver.

\section{Numerical results}
\label{hhbsec:gurobi}

We generate 87
random test problems\footnote{In \cite{hhbBK}, the authors compared
CycIP with a Swiss Army Knife. The Wenger Swiss Army Knife version
XXL, listed in the Guinness Book of World Records as the world's
most multi-functional penknife, contains 87 tools.} as in
\cite{hhbBK}. 
The size of each problem, $n$, satisfies
$341\leq n\leq 2735$. These problems are significantly larger
than those of \cite{hhbBK} because we wish to compare execution
time rather than number of iterations. 

Consider the following two measures of infeasibility:
\begin{equation}
(\forall x\in X)\quad
d_2(x) :=  \sqrt{\sum_{i=1}^{6} d_{C_i}^2(x)}
\quad\text{and}\quad
d_\infty(x) := \max_{i\in I}\|x-P_{C_i}x\|_\infty,
\end{equation}
where $\|x\|_\infty$ is the max-norm\footnote{Recall that if
$x=(\xi_1,\ldots,\xi_n)\in X$, then $\|x\|_\infty =
\max\{|\xi_1|,\ldots,|\xi_n|\}$.} of $x$.
Note that $d_2(x)=d_\infty(x)=0$ if and only if $x\in C$. 
Set $\varepsilon := 5\cdot 10^{-4}$, and let 
$(x_k)_{k\in\mathbb{N}}$ be a sequence generated by CycIP. 
We employ either 
$d_2(x_k)<\varepsilon$ or 
$d_\infty(x_k)<\varepsilon$ as stopping criterion. 

Let $\mathcal{P}$ be the set
of test problems, and let $\mathcal{A}$ be the set of algorithms. 
Let $(x_k^{(a,p)})_{k\in\mathbb{N}}$ be the 
sequence generated by algorithm $a\in\mathcal{A}$ applied to the problem
$p\in\mathcal{P}$. 
To compare the performance of the algorithms, 
we use \emph{performance profiles}\footnote{
For further information on performance profiles, 
we refer the reader to \cite{hhbDM}.}:
for every $a\in\mathcal{A}$ and 
for every $p\in\mathcal{P}$, we set
\begin{equation}
r_{a,p} :=
\frac{\tau_{a,p}}{\min\hhbmenge{\tau_{a',p}}{a'\in\mathcal{A}}} \geq 1,
\end{equation}
where $\tau_{a,p}\in\{1,2,\ldots,\tau_{\max}\}$ is the time that
$a$ requires to solve $p$ and $\tau_{\max}$ is the maximum time
allotted for all algorithms. 
If $r_{a,p} = 1$, then $a$ uses the least amount of time to
solve problem $p$. 
If $r_{a,p} > 1$, then $a$ requires $r_{a,p}$
times more time for $p$ than the algorithm that uses the least
amount of time for $p$.
For each algorithm $a\in\mathcal{A}$, we plot the function
\begin{equation}
\rho_a\colon \mathbb{R}_+\to[0,1]\colon \kappa\mapsto 
\frac{\hhbcard\hhbmenge{p\in\mathcal{P}}{\log_2(r_{a,p})\leq\kappa}}{\hhbcard
\mathcal{P}},
\end{equation}
where ``$\hhbcard$'' denotes the cardinality of a set. 
Thus, $\rho_a(\kappa)$ is the percentage of problems that algorithm
$a$ solves within factor $2^\kappa$ of the best algorithms. 
Therefore, an algorithm $a\in\mathcal{A}$ 
is ``fast'' if $\rho_a(\kappa)$ is large for $\kappa$ small;
and $a$ is ``robust'' if $\rho_a(\kappa)$ is large for $\kappa$ large. 

To compare CycIP with a linear programming solver, we model \eqref{hhbe:road}
as the constraints of a \emph{Linear Program (LP)}. As objective
function, we use $x\mapsto \|x-x_0\|_1$, where $\|x\|_1$ denotes
the $1$-norm\footnote{Recall that if $x=(\xi_1,\ldots,\xi_n)\in
X$, then $\|x\|_1 = |\xi_1|+\cdots+|\xi_n|$.} of $x$.
As LP solver, we use Gurobi 5.5.0, 
a state-of-the-art mathematical programming solver \cite{hhbGurobi}. 
CycIP was implemented with the C++ programming language.
We run the experiments on a Linux computer with a 2.4 GHz Intel\textsuperscript{\textregistered} Xeon\textsuperscript{\textregistered} E5620 CPU
and 24 GB of RAM. As time measurement, we use \emph{wall-clock
time}\footnote{To allow for a more fair comparison,
we included in wall-clock time only the time
required for running the solver's software itself (and not
the time for loading the problem data or for setting up the solver's 
parameters).}.
We limit the solving time to
$\tau_{\max} := 150$ seconds for each problem and algorithm. 

\begin{figure}[h!]
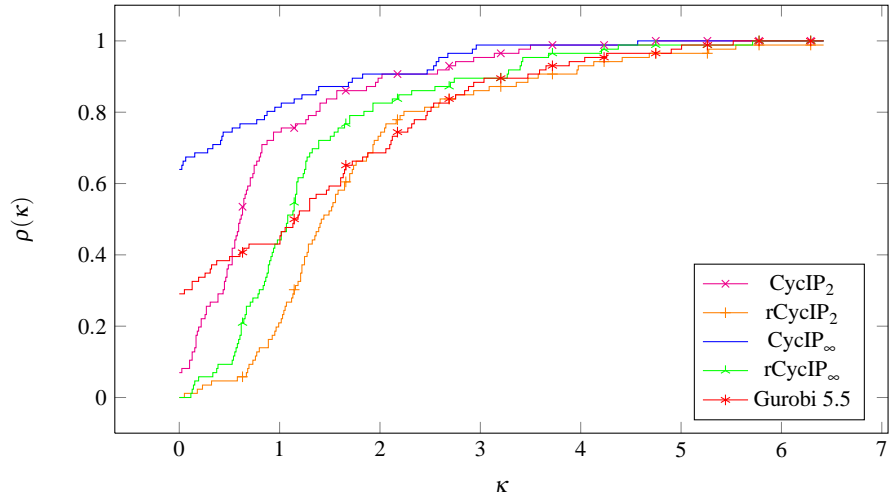

\centering

\caption{Performance profile for convex problems.\label{hhbfig:pprof_convex}}
\end{figure}

Figure~\ref{hhbfig:pprof_convex} shows the performance profile
for the convex case. 
Here, CycIP uses a cyclic control with period $6$
and the randomized \textbf{rCycIP} variant has quasicyclic control 
satisfying $(\forall k\in\mathbb{N})$
$\{\mathrm{i}(6k),\mathrm{i}(6k+1),\ldots,\mathrm{i}(6k+5)\} =
\{1,2,\ldots,6\}$, i.e., 
for every $k\in\mathbb{N}$, 
$(\mathrm{i}(6k),\mathrm{i}(6k+1),\ldots,\mathrm{i}(6k+5))$ is a
randomly generated permutation of $(1,2,\ldots,6)$. 
Depending on whether $d_2$ or $d_\infty$ was used as the 
infeasibility measure, we write 
$\mathrm{CycIP}_2$ and 
$\mathrm{CycIP}_\infty$, respectively, and similarly for rCycIP. 

\FloatBarrier
\FloatBarrier

\begin{figure}[h!]
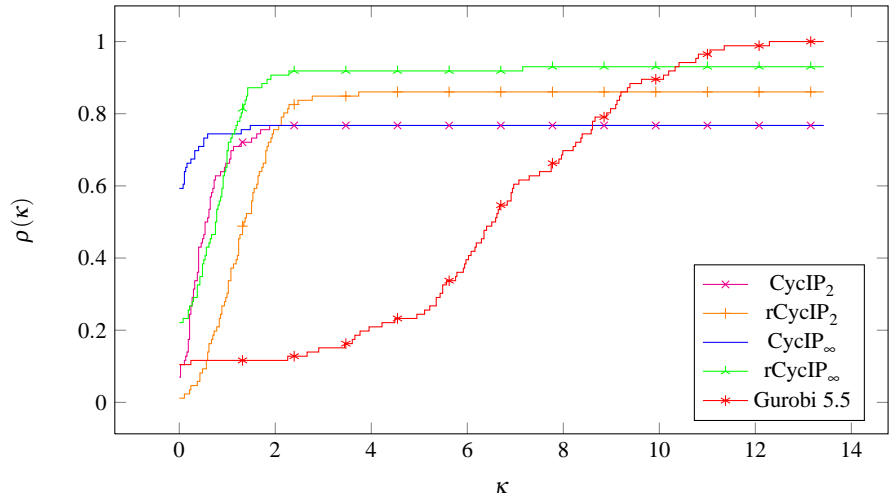

\centering
%
\caption{Performance profile for nonconvex problems.\label{hhbfig:pprof_nonconvex}}
\end{figure}

For the nonconvex case, shown in Figure~\ref{hhbfig:pprof_nonconvex}, we
included a minimum slope constraint as mentioned in 
Remark~\ref{hhbrem:nonconvex}.
For the LP solver, the resulting model becomes a \emph{Mixed Integer Linear Program}.

\FloatBarrier
\FloatBarrier

We infer from the figures that for convex problems,
$\mathrm{CycIP}_\infty$ solves the test problems quickly and
robustly. For nonconvex problems, $\mathrm{CycIP}_\infty$ is
still fast, but less robust than the slower randomized
variant $\mathrm{rCycP}_\infty$. 
Gurobi is the slowest --- but also the most robust --- algorithm. 

\section{Conclusion}

\label{hhbsec:conclusion}

In this work, we proved that the method of cyclic intrepid
projections converges to a feasible solution under quasicyclic
control and an interiority assumption.
Specialized to a problem arising in road design, this leads to
the first rigorous proof of convergence of CycIP.
Numerical results show that CycIP
is competitive compared to a commercial optimization solver,
especially in terms of speed. 
Randomization strategies increase robustness in case of nonconvex
problems for which there is no underlying convergence theory. 
Future work will focus on obtaining theoretical convergence results 
and on experimenting with other algorithms to
increase robustness in the nonconvex setting. 

\begin{acknowledgement}
The authors thank Dr.~Ramon Lawrence for the opportunity to run 
the numerical experiments on his server, and 
Scott Fazackerley and Wade Klaver for technical help.
HHB also thanks Dr.\ Masato Wakayama and the 
Institute of Mathematics for Industry, Kyushu University,
Fukuoka, Japan for their hospitality --- some of this research
benefited from the extremely stimulating environment during the 
``Math-for-Industry 2013'' forum.
HHB was partially supported by the Natural Sciences and
Engineering Research Council of Canada (Discovery Grant and
Accelerator Supplement) and
by the Canada Research Chair Program.
\end{acknowledgement}

% %%

%%%%%%%%%%%%%%%%%%%%%%%%%%%%%%%%%%%%%%%%%%%%%%%%%%%%%%%%%%%%%%%%%%%%%%

\end{document}